\documentstyle{article}
\title{A New Aspect of Representations of $U_q(\hat{sl_2})$\\
---Root of Unity Case}
\author{Xufeng Liu\\
Department of Mathematics, The University of Melbourne, Parkville, Victoria 3052,
Australia\\
Department of Mathematics, Peking University, Beijing 100871, P.R.China
\footnote{Permanent mailing address}}
\begin{document}
\begin{center}
\begin{flushleft}
{\bf A New Aspect of Representations of $U_q(\hat{sl_2})$\\
---Root of Unity Case}\\
\vspace{4mm}
Xufeng Liu\\
\vspace{3mm}
Department of Mathematics, The University of Melbourne\\ Parkville, Victoria 3052,
Australia\\
\vspace{3mm}
Department of Mathematics, Peking University\\ Beijing 100871, P.R.China
\footnote{Permanent mailing address}
\end{flushleft}
\end{center}
\vspace{20mm}
\begin{flushleft}
Abstract.The structure of the tensor product representation $V_{\lambda _{1}}(x)
\otimes V_{\lambda _{2}}(y)$ of $U_q(\hat{sl_2})$ is probed at roots of unity.
A polynomial  identity is derived as an outcome.
Also, new bases of $V_{\lambda _{1}}(x)\otimes V_{\lambda _{2}}(y)$
are established under certain conditions.
\end{flushleft}
\newpage
In this paper
we take $q$ to be a complex number instead of a formal variable,
unless otherwise pointed out explicitly.For any complex
numbers $q$ and $\lambda$ and any integers $r$ and $s$ we use the notations
\begin{eqnarray*}  [\lambda]_q&=&\frac{q^{\lambda}-q^{\lambda}}{q-q^{-1}},\\
 \left [r\right ]_{q}!&=&\prod_{s=1}^{r} [s]_{q},\\
\left [\begin{array}{c}
r\\s \end{array}\right ]&=&\frac{[r]_{q}!}{[s]_{q}![r-s]_{q}!},r\geq s.
\end{eqnarray*}
and we set $[0]_{q}!=1$.

\vspace{5mm} 
{\bf 1.Introduction and Some Basic Facts}\\
In recent years, the representation theory of quantum affine algebras has been
well developed. In a series of papersp[1,2,3], Chari and Pressley establish
the finite dimensional representation theory when q is not a root of unity
and in [4] they classify a certain category of irreducible finite dimensional
representations of the so-called restricted quantum affine algebras at roots
of unity. On the other hand, the case of non-restricted quantum affine
algebras has been dealt with by Kac and Beck[5]. In this paper, we will consider
some non-standard representations, whose meaning will become clear
 as we proceed, of the simplest quantum affine algebra
$U_q(\hat{sl_2})$, in the
non-generic case.Most unexpectedly we are able to derive some polynomial
identities from the representation theory.

Let us now  list some basic facts about $U_q(\hat{sl_2})$ and fix
notations.For details we refer the readers to Ref.[6].

Definition 1.1. The quantum affine algebra $U_q(\hat{sl_2})$ 
 is the associative algebra over $C$ with generators $e_{i},f_{i},K_{i}$ and
$K_{i}^{-1}$($i=0,1$) and the following relations:
\begin{eqnarray*}
&&K_{i}K_{i}^{-1} = K_{i}^{-1}K_{i}=1,\\
&&K_0K_1 = K_1K_0,\\
&&K_{i}e_{i}K_{i}^{-1}=q^2e_{i},K_{i}f_{i}K_{i}^{-1}=q^{-2}f_{i},\\
&&K_{i}e_{j}K_{i}^{-1}=q^{-2}e_{j},K_{i}f_{j}K_{i}^{-1}=q^2f_{j},i\ne j, \\
&&\left [e_{i},f_{i}\right ]=\frac{K_{i}-K_{i}^{-1}}{q-q^{-1}},\\
&&\left [e_0,f_1\right ]=\left [e_1,f_0\right ]=0,
\end{eqnarray*}
\begin{eqnarray*}
e_{i}^3e_{j}-[3]_qe_{i}^2e_{j}e_{i}+[3]_qe_{i}e_{j}e_{i}^2-e_{j}e_{i}^3&=&0,\\
f_{i}^3f_{j}-[3]_qf_{i}^2f_{j}f_{i}+[3]_qf_{i}f_{j}f_{i}^2-f_{j}f_{i}^3&=&0,(i\ne j).
\end{eqnarray*}
Moreover, $U_q(\hat{sl_2})$ is a Hopf algebra over $C$ with the comultiplication
\begin{eqnarray*}
&&\triangle (e_{i})=e_{i}\otimes K_{i}+1\otimes e_{i},\\
&&\triangle (f_{i})=f_{i}\otimes 1+K_{i}^{-1}\otimes f_{i},\\
&&\triangle (K_{i})=K_{i}\otimes K_{i},\triangle (K_{i}^{-1})=K_{i}^{-1}\otimes K_{i}^{-1},\\
\end{eqnarray*}
and the antipode
\begin{eqnarray*}
&&S(K_{i})=K_{i}^{-1},S(K_{i}^{-1})=K_{i},\\
&&S(e_{i})=-e_{i}K_{i}^{-1},S(f_{i})=-K_{i}f_{i}.
\end{eqnarray*}

Definition 2.2. The quantum algebra $U_q(sl_2)$ is the associative algebra over
$C$ with generators $e,f,K$ and $K^{-1}$ and the following relations:
\begin{eqnarray*}
&&KK^{-1}=K^{-1}K=1,\\
&&KeK^{-1}=q^2e,KfK^{-1}=q^{-2}f,\\
&&\left [e,f\right ]=\frac{K-K^{-1}}{q-q^{-1}}.
\end{eqnarray*}
It is a Hopf algebra over $C$ with the following comultiplication $\triangle$
and antipode $S$:
\begin{eqnarray*}
&&\triangle e=e\otimes K+1\otimes e,\\
&&\triangle f=f\otimes 1+K^{-1}\otimes f,\\
&&\triangle K=K\otimes K,\triangle K^{-1}=K^{-1}\otimes K^{-1},\\
&&S(K)=K^{-1},S(K^{-1})=K,S(e)=-eK^{-1},S(f)=-Kf.
\end{eqnarray*}
 There is an associative algebra homomorphism, known as evaluation map, from
$U_q(\hat{sl_2})$ to $U_q(sl_2)$[7].

Definition 1.3.For any $x\in C\backslash\{0\}$, the evaluation map $ev_x$ 
 from $U_q(\hat{sl_2})$ to $U_q(sl_2)$ is the associative algebra homomorphism
such that
\begin{eqnarray*}
&&ev_x(e_0)=q^{-1}xf, ev_x(e_1)=e,\\
&&ev_x(f_0)=qx^{-1}e,ev_x(f_1)=f,\\
&&ev_x(K_0)=K^{-1},ev_x(K_1)=K.
\end{eqnarray*}
It is clear that modules of $U_q(\hat{sl_2})$ can be obtained by pulling back
modules of $U_q(sl_2)$ by the homomorphism $ev_x$.We denote by $V(x)$ the pull
back module of $U_q(\hat{sl_2})$ of a module $V$ of $U_q(sl_2)$ by the evaluation
map $ev_x$.

Definition 1.4. Let $V$ and $W$ be two modules of $U_q(sl_2)$.The $U_q(\hat{sl_2})$ module
$W(y)\otimes V(x)$ is the vector space $W\otimes V$ with the module structure
defined through the following action
\[g(w\otimes v)\stackrel{\triangle}{=}(ev_y\otimes ev_x)\triangle g(w\otimes
v),\forall g\in U_q(\hat{sl_2}).\]

{\bf Structure of $V_{\lambda _{1}}(x)\otimes V_{\lambda _{2}}(y)$}

This section is devoted to the study of tensor product
representations of $U_q(\hat{sl_2})$ arising from non-standard representations
of $U_q(sl_2)$ when $q$ is a root of unity. From now on we take $q$ to be a
root of unity. Let $p$ be the smallest positive integer such that
$q^{p}=1$.For simplicity, we assume that $p$ is odd.

It is best known that for any complex number $\lambda$ there is a $p$-dimensional
$U_q(sl_2)$ module $V=span\{v_{i}|i=0,1,2,\cdots p-1\}$ such that
\begin{eqnarray*}
&&Kv_{i}=q^{\lambda-2i}v_{i},K^{-1}v_{i}=q^{-(\lambda -2i)}v_{i},i=0,1,\cdots p-1,\\
&&fv_{i}=[i+1]v_{i+1},fv_{p-1}=0,i\ne p-1,\\
&&ev_{i}=[\lambda +1-i]v_{i-1},ev_{0}=0,i\ne 0.
\end{eqnarray*}
Obviously,if
$$
[\lambda +1-i]\ne 0,i=1,2,\cdots p-1
$$
this representation is irreducible. We denote it by $V_{\lambda}$,and call
the basis satisfying the above equations normal basis.

Let $\lambda _{1}$ and $\lambda _{2}$ be two complex numbers.
 We will investigate the struture of $V_{\lambda _{1}}(x)\otimes
V_{\lambda _{2}}(y)$.
For convenience we introduce the following notions concerning this
representation.

Condition 1.\begin{eqnarray*}
&&[\lambda _{j}+1-i]\ne 0,i=0,1,2,\cdots p,j=1,2,\\
&&\left [\lambda _{1}+\lambda _{2}-2j+1-i\right ]\ne 0,i=0,1,2,\cdots p,
j=0,1,\cdots p-1.
\end{eqnarray*}

Condition 2.$y/x\neq q^{\lambda _{1}+\lambda _{2}-2(l-1)},l=1,2,\cdots,p-1.$

Condition 3.$y/x\neq q^{-(\lambda _{1}+\lambda _{2}-2(l-1))},l=1,2,\cdots,p-1.$

We will need the following elementary result.

Lemma 2.1.Let $V$ be a $U_{q}(sl_2)$ module. If $\omega$ is a highest
weight vector of weight $\lambda$ and $\omega \prime$ is a lowest vector of
weight $\lambda \prime$,we have
\begin{eqnarray*}
&&ef^{n}\omega =[n]_{q}[\lambda -n+1]_{q}f^{n-1}\omega\\
&&fe^{n}\omega \prime =-[n]_{q}[\lambda \prime +n-1]_{q}e^{n-1}\omega \prime.
\end{eqnarray*}

We denote by $\{v_{i}|i=0,1,\cdots,p-1\}$ and $\{w_{i}|i=0,1,\cdots,p-1\}$
the normal bases of $V_{\lambda _{1}}$ and $V_{\lambda _{2}}$ respectively,
and for each $l\in \{0,1,\cdots,p-1\}$ define
\begin{eqnarray*}
&&\Omega _{l}=\sum_{i=0}^{l} c_{i,l-i}v_{i}\otimes w_{l-i},\\
&&\Phi _{l}=\sum_{i=l}{p-1} d_{i,l+p-1-i}v_{i}\otimes w_{l+p-1-i},0\leq l\leq
p-1\end{eqnarray*}
where
\begin{eqnarray*}
&&c_{i,l-i}=(-1)^{i}q^{i(2l-i-\lambda _{2}-1)}\prod_{j=0}^{i}
\frac{[\lambda _{2}-l+j]_{q}}{[\lambda _{1}-j+1]_{q}},\\
&&d_{i,l+p-1-i}=(-1)^{i-l}q^{(i-l)(\lambda _{1}-i-l-1)}\prod_{j=l+1}{i}
\frac{[j]_{q}}{[l+p-j]_{q}},i\geq l+1,\\
&&d_{l,p-1}=[l]_{q},\end{eqnarray*}
when they make sense.

Lemma 2.2.Under Condition 1 $\Omega _{l}(l=0,1,\cdots,p-1)$ are the only
highest weight vectors of $U_q(sl_2)$ (up to scalar multiples) in
$V_{\lambda _{1}}\otimes V_{\lambda _{2}}$ and $\Phi _{l} $ are the only
lowest vectors.

Proof.Condition 1 guarantees their well definedness. The fact that they
are highest weight and lowest weight vectors respectively and the
uniqueness can be easily established by direct calculation.

Lemma 2.3.Presuppose Condition 1. Let k be the smallest positive integer
such that $f^{k}\Omega _{l}=0$ but $f^{k-1}\Omega _{l}\ne 0$.We have
$k=p$.

Proof.That $k$ cannot be smaller than $p$ follows from Lemma 2.1.On the other
hand,if $k$ were larger than $p$, $f^{p}\Omega _{l}$ would be a highest
weight vector. Again from Lemma 2.1 this would imply $\lambda _{1}+
\lambda _{2}-2l-2p-2(p-1)$ is a weight.But this contradicts the fact that
the lowest weight in $V_{\lambda _{1}}\otimes V_{\lambda _{2}}$ is
$\lambda _{1}+\lambda _{2}-4(p-1)$.

Proposition 2.1.If Condition 1 is satisfied 
 $V_{\lambda _{1}}(x)\otimes V_{\lambda _{2}}(y)$ ,as $U_q(sl_2)$ module,
 has the following decomposition:
$$
V_{\lambda _{1}}(x)\otimes V_{\lambda _{2}}(y)=
\sum_{j=0}^{p-1}\oplus V_{\lambda _{1}+\lambda _{2}-2j}
$$

Proof.From Lemma 2.3 the submodule $U_q(sl_2)\Omega _{l}$ is spanned
by $\{f^{i}\Omega _{l}|i=0,1,\cdots,p-1\}$ and hence has dimension $p$.
Lemma 2.1 says this submodule is isomorphic to $V_{\lambda _{1}
+\lambda _{2}-2l}$.On the other hand, one can easily prove, by means of
the raising operator $e$, that vectors from different $U_{q}(sl_2)$
are linearly independent.Now the proposition follows as a result of
simple dimension counting.

From now on we will presuppose Condition 1.So we can use Proposition 2.1
freely.

According to Lemmas 2.2 and 2.3 for each $l$ $f^{p-1}\Omega _{l}$ must be
a scalar multiple of $\Phi _{l}$.We will denote this scalar by $\alpha _{l}$:
$f^{p-1}\Omega _{l}=\alpha _{l}\Phi _{l}.$

To establish the main theorem of this section we need more lemmas. The
following two lemmas can be verified by direct calculation.

Lemma 2.4.\begin{eqnarray*}
f_{0}\Omega _{l}& = &q[\lambda _{2}-l]_{q}(-x^{-1}q^{-\lambda _{2}+
2l-2}+y^{-1}q^{\lambda _{1}})\Omega _{l-1},\\
e_{0}\Phi _{l} & = & q^{-3}[l]_{q}(xq^{-\lambda _{2}}-yq^{\lambda _{1}-2l})
\Phi _{l+1}.
\end{eqnarray*}

Lemma 2.5.\begin{eqnarray*}
e^{2}e_{0}\Omega _{l} & = & q^{-1}[2]_{q}[\lambda _{2}-l]_{q}
(xq^{\lambda _{1}}-yq^{-\lambda_{2}+2l-2})\Omega _{l-1},\\
f^{2}f_{0}\Phi _{l} & = & [2]_{q}[l]_{q}(-x^{-1}q^{\lambda _{1}-2l-1}
+y^{-1}q^{-\lambda _{2}-1})\Phi _{l+1}.
\end{eqnarray*}

Here in stating these lemmas we have tacitly adopted the convention
$\Omega _{1}=\Phi _{p}=0$.Thanks to Lemma 2.5 and Proposition 2.1
 we can write down
\begin{eqnarray*}
e_{0}\Omega _{l} & = & \beta _{l,l-1}f^{2}\Omega _{l-1}+\beta _{l,l}f\Omega {l}
+\beta _{l,l+1}\Omega _{l+1},\\
f_{0}\Phi _{l} & = & \gamma _{l,l+1}e^{2}\Phi _{l+1}+\gamma _{l,l}
e\Phi _{l}+\gamma _{l,l-1}\Phi _{l-1}.
\end{eqnarray*}
where the coefficients are some complex numbers.Before proceeding along
let us introduce another notation:for an arbitrary complex number $z$
$$\overline{[z]_{q}}=\frac{q^{z}-q^{-z}}{q-q^{-1}}.$$
Now for two complex numbers $z_{1}$ ,$z_{2}$ and an integer $l$ define
$$ A_{l}(z_{1},z_{2})=
 \overline{[z_{1}+1]_{q}}
\overline{[z_{1}+z_{2}-2l+1]_{q}}-\overline{[1]_{q}}
\overline{[z_{2}+1]_{q}}$$

Lemma 2.6.\begin{eqnarray*}
\beta _{l,l-1} & = & \frac{q^{-1}[\lambda _{2}-l]_{q}(xq^{\lambda _{1}}
-yq^{-\lambda _{2}+2l-2})}{[\lambda _{1}+\lambda _{2}-2l+1]_{q}
[\lambda _{1}+\lambda _{2}-2l+2]_{q}},\\
\beta _{l,l} & = & \frac{xq^{-1}A_{l}(\lambda _{1},\lambda _{2})
+yq^{-1}A_{l}(\lambda _{2},\lambda _{1})}{[\lambda _{1}+\lambda _{2}-2l]_{q}
[\lambda _{1}+\lambda _{2}-2l+2]_{q}},\\
\beta _{l,l+1} & = & -\frac{[l+1]_{q}[\lambda _{2}-l]_{q}[\lambda _{1}-l]_q
[\lambda _{1}+\lambda _{2}-l+1]_{q}}{[\lambda _{1}+\lambda _{2}-2l]_{q}
[\lambda _{1}+\lambda _{2}-2l+1]_{q}}(xq^{-\lambda _{1}}-yq^{\lambda _{2}-2l)}
q^{-1},\\
\beta _{0,0} & = & \frac{q^{-1}}{[\lambda _{1}+\lambda _{2}]_{q}}
(x[\lambda _{1}]_{q}+y[\lambda _{2}]_{q}),\\
\beta _{0,1} & = & \frac{q^{-1}[\lambda _{1}]_{q}[\lambda _{2}]_{q}}
{[\lambda _{1}+\lambda _{2}]_{q}[\lambda _{2}-1]_{q}}(-xq^{-\lambda _{1}}
+yq^{\lambda _{2}}).
\end{eqnarray*}

Proof.$\beta _{l,l-1}$ can be determined easily by using the first
equation of Lemma 2.5.$\beta _{l,l}$ can be derived in a similar way.
Finally $\beta _{l,l+1}$ can be obtained from them by direct calculation.

Similarly we can prove the following

Lemma 2.7.\begin{eqnarray*}
\gamma _{l,l+1} & = & \frac{q[l]_{q}(-x^{-1}q^{\lambda _{1}-2l-2}
+y^{-1}q^{-\lambda _{2}})}{[\lambda _{1}+\lambda _{2}-2l]_{q}
[\lambda _{1}+\lambda _{2}-2l+1]_{q}},\\
\gamma _{l,l} & = & \frac{xqA_{l}(\lambda _{1},\lambda _{2})
+yqA_{l}(\lambda _{2},\lambda _{1})}{[\lambda _{1}+\lambda _{2}-2l]_{q}
[\lambda _{1}+\lambda _{2}-2l+2]_{q}},\\
\gamma _{l,l-1} & = & -\frac{[l]_{q}[\lambda _{2}-l+1]_{q}
[\lambda _{1}-l+1]_q
[\lambda _{1}+\lambda _{2}-l+2]_{q}}{[l-1]_{q}[\lambda _{1}+\lambda _{2}-2l+1]_{q}
[\lambda _{1}+\lambda _{2}-2l+2]_{q}}(-x^{-1}q^{-\lambda _{1}+2l}
+y^{-1}q^{\lambda _{2}+2})q,\\
\gamma _{p-1,p-1} & = & \frac{q}{[\lambda _{1}+\lambda _{2}+4]_{q}}
(x^{-1}[\lambda _{1}+2]_{q}+y^{-1}[\lambda _{2}+2]_{q}),\\
\gamma _{p-1,p-2} & = & \frac{q^{3}[\lambda _{1}+2]_{q}[\lambda _{2}+2]_{q}}
{[\lambda _{1}+\lambda _{2}+3]_{q}[2]_{q}}(x^{-1}q^{-\lambda _{1}-4}
-y^{-1}q^{\lambda _{2}}).
\end{eqnarray*}

Proposition 2.2.For each $1\leq l\leq p-1$when $y/x=q^{\lambda _{1}+
\lambda _{2}-2(l-1)}$, $\sum_{j=l}^{p-1}\oplus V_{\lambda _{1}+\lambda _{2}
-2j}$ is a $U_q(\hat{sl_2})$ submodule of $V_{\lambda _{1}}(x)\otimes
V_{\lambda _{2}}(y)$;when$y/x=q^{-(\lambda _{1}+
\lambda _{2}-2(l-1))}$, $\sum_{j=0}^{l}\oplus V_{\lambda _{1}+\lambda _{2}
-2j}$ is a $U_q(\hat{sl_2})$ submodule of $V_{\lambda _{1}}(x)\otimes
V_{\lambda _{2}}(y)$.

Proof. We notice that $V_{\lambda _{1}+\lambda _{2}}$ can be generated by
$\Omega _{l}$ or $\Phi _{l}$. With this fact in mind one easily sees that
this proposition follows directly from Lemmas 2.4,2.6 and 2.7.

Similarly one can prove the following

Proposition 2.3.If Condition 2 is satisfied by
$V_{\lambda _{1}}(x)\otimes V_{\lambda _{2}}(y)$ it has no submodule
not containing the highest weight vector $\Omega _{0}$ and if Condition 3
is satisfied it has no proper submodule containing $\Omega _{0}$.

We are now finally in a position to prove the following main result of
this section.

Theorem 2.1.   $V_{\lambda _{1}}(x)\otimes
V_{\lambda _{2}}(y)$ is irreducible as  $U_q(\hat{sl_2})$ module if and only
if Condition 2 and Condition 3 are satisfied.

Proof.The "if" part is an easy consequence of the observation that under
the condition of this theorem,$V_{\lambda _{1}}(x)\otimes V_{\lambda _{2}}
(y)$ has no proper submodule, owing to 
 Proposition 2.3. The "only if" part is just another version of
Proposition 2.2.

{\bf 3.An Identity}

This section is devoted to establish some identities.The trick is simple:
we determine the coefficient $\alpha _{l}$ in two different ways and then
equate the results.

For any positive integer $k$ we have
$$
\triangle f^{k}=\sum_{j=0}{k} q^{-j(k-j)}\left [\begin{array}{c}
k\\j\end{array}\right ]_{q}K^{-j}f^{k-j}\otimes f^{j}.
$$
Using this formula, after some elementary calculation, we get the following
coefficient of the term $v_{p-l}\otimes w_{l}$ in $f^{p-1}\Omega _{l}$:
$$
\sum_{j=0}^{l}q^{-j(\lambda _{1}-j+1)}\frac{([p-1]_{q}!)^{2}[l]_{q}!}
{([j]_{q}!)^{2}[l-j]_{q}![p-1-j]_{q}!}c_{j,l-j}.
$$
On the other hand, the coefficient of the same term in $\Phi _{l} $ is
$$
d_{p-1,l}=(-1)^{l}[l]_{q}q^{l(l+1)}q^{\lambda _{1}(p-1-l)}.
$$
It follows that
$$
\alpha _{l}=(-1)^{l}q^{-l(l+1)}q^{-\lambda _{1}(p-1-l)}[l-1]_{q}!
\sum_{j=0}^{l}q^{-j(\lambda _{1}-j+1}\frac{([p-1]_{q}!)^{2}}
{([j]_{q}!)^{2}[l-j]_{q}![p-1-j]_{q}!}c_{j,l-j}.
$$

Now let us turn to another way of determining $\alpha _{l}$.
Since for each $l\in \{0,1,\cdots,p-1\}$
$f^{p}\Omega _{l}=0$ we have the equation
$$e_{0}f^{p-1}\Omega _{l}=\beta _{l,l+1}f^{p-1}\Omega _{l+1}.$$
It then follows from the second equation of Lemma 2.4 and Lemma 2.6 that
\begin{eqnarray*}
\alpha _{l+1}/\alpha _{l} & = & q^{-3}[l]_{q}(xq^{-\lambda _{2}}
-yq^{\lambda _{1}-2l}) \\
& = & -q^{\lambda _{1}-\lambda _{2}-2}\frac{[l]_{q}[\lambda _{2}-l-1]_{q}
[\lambda _{1}+\lambda _{2}-2l]_{q}[\lambda _{1}+\lambda _{2}-2l+1]_{q}}
{[l+1]_{q}[\lambda _{2}-l]_{q}[\lambda _{1}-l]_{q}[\lambda _{1}+
\lambda _{2}-l+1]_{q}},\\
\frac{\alpha _{1}}{\alpha _{0}} & = & (-1)^{p}q^{\lambda _{1}-\lambda _{2}
-2}\frac{[\lambda _{2}-1]_{q}[\lambda _{1}+\lambda _{2}]_{q}}
{[\lambda _{2}]_{q}[\lambda _{1}]_{q}}.
\end{eqnarray*}
Now $\alpha _{l}$ can be determined as follows.
\begin{eqnarray*}
\alpha _{l} & = & \frac{\alpha _{l}}{\alpha _{0}}=\alpha _{0}
\prod_{j=1}^{l} \frac{\alpha _{j}}{\alpha _{j-1}}\\
& = & (-1)^{p+l}q^{-(p-1)\lambda _{1}}q^{(\lambda _{1}-\lambda _{2}-2)l}
\frac{[l-1]_{q}![\lambda _{2}-l]_{q}[p-1]_{q}!}{[l]_{q}![\lambda _{1}+1]_{q}
[\lambda _{1}-l+1]_{q}}\times\\
& & \frac{\prod_{j=0}^{2(l-1)} [\lambda _{1}+\lambda _{2}-2l+j+2]_{q}}
{\prod_{j=1}{l-1} [\lambda _{1}+\lambda _{2}-l+j+1]_{q}[\lambda _{1}
-l+j+1]_{q}}.
\end{eqnarray*}
Here we have used
$$\alpha _{0}=-(-1)^{p}q^{-(p-1)\lambda _{1}}\frac{[p-1]_{q}!
[\lambda _{2}]_{q}}{[\lambda _{1}+1]_{q}}.$$

Equating this result and the previous one and simplifying the equation,
we arrive at the following conclusion: for each $l\leq p-1$
\begin{eqnarray*}
&&q^{l\lambda _{2}}\sum_{j=0}^{l}q^{-j(\lambda _{1}+\lambda _{2})}q^{2jl}
q^{-2j}\left [\begin{array}{c}l\\j\end{array}\right ]_{q}
\prod_{i=1}^{j} [\lambda _{2}-l+i]_{q}\prod_{i=1}^{l-j} [\lambda _{1}
-l+i]_{q}\\
&&=q^{l(l-1)}\prod_{i=0}^{l-1} [\lambda _{1}+\lambda _{2}-2l+j+2]_{q}.
\end{eqnarray*}
Here we set
$$\prod_{i=1}^{0} [\lambda _{2}-l+i]_{q}=\prod_{i=1}^{0} [\lambda _{1}
-l+i]_{q}=1.$$

Examining the process we know that this equation is subject to the restraints
$q^{p}=1$,$l\leq p-1$ and Condition 1.But it is only specious.Indeed the
restraints can be lifted easily by the following argument.Regarding both
sides as continuous functions of $\lambda _{1}$ and $\lambda _{2}$ we can
legally discard Condition 1. Especially the equation is true when
$\lambda _{1}$ and $\lambda _{2}$ are integers.Now by specializing
$\lambda _{1}$ and $\lambda _{2}$ to some integers we think of the equation
as a polynomial equation in indeterminate $q$. Then if a complex number $z$
satisfies the condition:
$ z^{p}=1,p$ odd and $p\geq l$ it is a root of the equation. There are
obviously infinitely many such complex numbers. Therefore,the equation
must be a polynomial identity. To sum up we have proved the following

Theorem 3.1.Let $q,z_{1}$ and $z_{2}$ be arbitrary complex numbers.Then
for each positive number $l$
\begin{eqnarray*}
&&q^{lz_{2}}\sum_{j=0}^{l}q^{-j(z_{1}+z_{2})}q^{2jl}
q^{-2j}\left [\begin{array}{c}l\\j\end{array}\right ]_{q}
\prod_{i=1}^{j} [z_{2}-l+i]_{q}\prod_{i=1}{l-j} [z_{1}
-l+i]_{q}\\
&&=q^{l(l-1)}\prod_{i=0}{l-1} [z_{1}+z_{2}-2l+j+2]_{q}.
\end{eqnarray*}
Corollary 1.Keep the notations.We have
\begin{eqnarray*}
&&q^{lz_{2}}\sum_{j=0}^{l}q^{-j(z_{1}+z_{2})}q^{2jl}
q^{-2j}\left [\begin{array}{c}l\\j\end{array}\right ]_{q}
\prod_{i=1}^{j} [z_{2}-l+i]_{q}\prod_{i=1}^{l-j} [z_{1}
-l+i]_{q}\\
&&=q^{-lz_{2}}q^{2l(l+1)}\sum_{j=0}^{l}q^{j(z_{1}+z_{2})}q^{-2jl}
q^{2j}\left [\begin{array}{c}l\\j\end{array}\right ]_{q}
\prod_{i=1}^{j} [z_{2}-l+i]_{q}\prod_{i=1}^{l-j} [z_{1}
-l+i]_{q}.
\end{eqnarray*}

Proof. The right hand side of the above equation is symmetric with respect
to $z_{1}$ and $z_{2}$, so should the left hand side be.The corollary
follows from this observation directly.

If we take $z_{1}$ and $z_{2}$ to be natural numbers larger than or equal
to $l$ Theorem 3.1 can be rewritten in the following form.

Corollary 2.Let $m,n$ and $l$ be natural numbers.Assume $m,n\geq l$.Then
$$
\sum_{j=0}^{l} \frac{[m-l+j]_{q}![n-j]_{q}!}{[j]_{q}![l-j]_{q}!}
=\frac{[m-l]_{q}![n-l]_{q}![m+n-l+1]_{q}!}{[l]_{q}![m+n-2l+1]_{q}!}.$$

{\bf 4.Isomorphism Theorem}

In this section we will prove an isomorphism theorem between
$V_{\lambda _{1}}(x)\otimes V_{\lambda _{2}}(y)$ and
$V_{\lambda _{2}}(y)\otimes V_{\lambda _{1}}(x)$.To simplify the subsequent
discussion let us make a notation convention:if $s$ is a symbol standing
for something concerning $V_{\lambda _{1}}(x)\otimes V_{\lambda _{2}}(y)$ 
$s'$ will be used to denote its counterpart concerning
$V_{\lambda _{2}}(y)\otimes V_{\lambda _{1}}(x)$. For example we denote by
$\Omega _{l}\prime$ the counterpart of $\Omega _{l}$.According to this convention,
the expression of $s'$ can be obtained from that of $s$ by interchanging
$\lambda _{1}$ and $\lambda _{2}$,$x$ and $y$ and $v_{i}$ and $w_{i}$
simultaneously.

From Proposition 2.1  there is a $U_{q}(sl_2)$
module isomorphism between
$V_{\lambda _{1}}(x)\otimes V_{\lambda _{2}}(y)$ and
$V_{\lambda _{2}}(y)\otimes V_{\lambda _{1}}(x)$ such that $\Omega _{l}$
is sent to its counterpart, namely,$\Omega _{l}\prime$. We denote this
isomorphism by $I$, and denote by $P_{l}$ the projection from
$V_{\lambda _{1}}(x)\otimes V_{\lambda _{2}}(y)$ onto its subspace
$V_{\lambda _{1}+\lambda _{2}-2l}$.Now define
$$\check{R}=\sum_{l=0}^{p-1} c_{l}I\cdot P_{l}$$
where
$$
c_{l}=\prod_{j=0}^{l} \frac{[\lambda _{2}-j]_{q}(xq^{\lambda _{1}}-
yq^{-\lambda _{2}+2j-2})}{[\lambda _{1}-j]_{q}(-xq^{-\lambda _{1}+2j-2}
+yq^{\lambda _{2}})}.
$$
It is evident that if Conditions 2 and 3 are all satisfied $\check{R}$
is a well defined nonzero mapping from
$V_{\lambda _{1}}(x)\otimes V_{\lambda _{2}}(y)$ to
$V_{\lambda _{2}}(y)\otimes V_{\lambda _{1}}(x)$ .Actually it turns out
to be an isomorphism under these conditions.

Let us first prove that it is an intertwiner.To establish the equation
$\check{R}e_{0}=e_{0}\check{R}$ one need only to check the equation
$$
\check{R}e_{0}\Omega _{l}=e_{0}\check{R}\Omega _{l}$$
because $check{R}$ is a $U_{q}(sl_2)$ module homomorphism.This equation
is equivalent to the following equations:
\begin{eqnarray*}
&&c_{l-1}\beta _{l,l-1}=c_{l}\beta _{l,l-1}\prime,\beta _{l,l}=\beta _{l,l}
\prime,\\
&&c_{l+1}\beta _{l,l+1}=c_{l}\beta _{l,l+1}\prime,\end{eqnarray*}
which can be verified directly.

Next we consider the equation
$\check{R}f_{0}=f_{0}\check{R}$.This time it is enough to verify
$\check{R}f_{0}\Phi _{l}=f_{0}\check{R}\Phi _{l}$,which boils down to
\begin{eqnarray*}
&&\frac{c_{l+1}}{c_{l}}=\frac{\alpha _{l}\prime}{\alpha _{l+1}\prime}
\frac{\alpha _{l+1}}{\alpha _{l}}\frac{\gamma _{l,l+1}\prime}
{\gamma _{l,l+1}},\gamma _{l,l}=\gamma _{l,l}\prime,\\
&&\frac{c_{l}}{c_{l-1}}=\frac{\alpha _{l-1}\prime}{\alpha _{l}\prime}
\frac{\alpha _{l}}{\alpha _{l-1}}\frac{\gamma _{l,l-1}}
{\gamma _{l,l-1}\prime}.\end{eqnarray*}
One can check these equations without any difficulty.

We have proved that $\check{R}$ is a $U_{q}(\hat{sl_2})$ module homomorphism
In fact it is an isomorphism when Conditions 2 and 3 are satisfied.
This is because under these conditions both
$V_{\lambda _{1}}(x)\otimes V_{\lambda _{2}}(y)$ and
$V_{\lambda _{2}}(y)\otimes V_{\lambda _{1}}(x)$ are irreducible.Thus
the "if" part of the following theorem is established.

Theorem 4.1.The $U_{q}(\hat{sl_2})$ modules
$V_{\lambda _{1}}(x)\otimes V_{\lambda _{2}}(y)$ and
$V_{\lambda _{2}}(y)\otimes V_{\lambda _{1}}(x)$ are isomorphic if and only
if Conditions 2 and 3 are satisfied.

Proof. Only the "only if" part remains to be proved.It follows from Condition
1 that if Condition 2 is violated by 
$V_{\lambda _{1}}(x)\otimes V_{\lambda _{2}}(y)$ it is satisfied by
$V_{\lambda _{2}}(y)\otimes V_{\lambda _{1}}(x)$ .Thus according
to Proposition 2.3  in that case they cannot be isomorphic. A similar
argument proves that if Condition 3 is violated they cannot be
isomorphic either.This completes the proof of the theorem.

Remark. It is possible to prove this theorem via a detour through
$U_{q}^{res}(\hat{sl_2})$[4].Essentially it then would follow from the
triangular decomposition theorem of  $U_{q}^{res}(\hat{sl_2})$ after some
arguments.But there are many technicalities to be dealt with.The
proof presented here is desirable because of its directness.

{\bf 5.A New Basis}

In this section we will establish a new basis of  $V_{\lambda _{1}}(x)
\otimes V_{\lambda _{2}}(y)$
under some condition. 

Let $j,l$ be two non-negative integers and $j, l\leq p-1$.We introduce the
notation
$$
\phi _{l,j}=e_{0}^{l-j}f^{j}\Omega _{0}
$$
and define the sets
$$
\Delta _{l}=\left \{\phi _{l,j} |j=0,1,\cdots ,l\right \},
l=0,1,\cdots,p-1. $$

Obviously, elements from different $\Delta _{l}$ are linearly independent
as they belong to different weight spaces. We will prove for each
$l\in \{0,1,\cdots,p-1\}$ $\Delta _{l}$ is a linearly independent set under
some condition. To this end, we will calculate explicitly the determinant
of the coefficient matrix of $\Delta _{l}$ with respect to the linearly
independent set $\{v_{i}\otimes w_{j}|i+j=l\}$.

Let
$$
\phi _{l,j}=\sum_{i=0}^{l} \gamma _{l,j}^{i,l-i}v_{i}\otimes w_{l-i}.
$$
We denote by $(\Delta _{l})$ the $l+1$ by $l+1$ coefficient matrix
$(\gamma _{l,j}^{i,l-i})$ whose row is marked by $j$ and column by
$(i,l-i)$ and denote by $|\Delta _{l}|$ the corresponding determinant.

Lemma 5.1.$$
\left |\Delta _{l+1}\right |=\frac{[\lambda _{2}]_{q}}
{[\lambda _{1}+1]_{q}}q^{-\frac{(l+1)(l+2)}
{2}}\prod_{j=1}^{l+1} [j]_{q}!\prod_{j=0}^{l}
(y-xq^{-\lambda _{1}-\lambda _{2}+2j})^{l+1-j}.$$

Proof.Repeat the proof verbatim of Proposition 3.1 in [8] with $m$ replaced
by $\lambda _{1}$ and $n$ replaced by $\lambda _{2}$.

Let us proceed to prepare another lemma. For each $l\in \{0,1,\cdots,p-1\}$
we define
$$\Lambda _{l}=\left \{e_{0}^{p-1-i}f^{i+l}\Omega _{0}|i=0,1,\cdots,p-1-l
\right \}$$
.For convenience we denote
$e_{0}^{p-1-i}f^{i+l}\Omega _{0}$ by $\varphi _{l,i}$. Then we have
$\varphi _{l+1,i}=f\varphi _{l,i}$.

Lemma 5.2. If Condition 3 is satisfied, then for
each $l\in \{0,1,\cdots,p-1\}$
both $\Delta _{l}$ and $\Lambda _{l}$ are linearly independent sets.

Proof.That $\Delta _{l}$ is a linearly independent set
follows directly from Lemma 5.1.
To prove that $\Lambda _{l}$ is also such a set we use induction method.
By definition $\Lambda _{0}=\Delta _{p-1}$.So when $l=0$ there is nothing
to be proved.Suppose $\Lambda _{j}$ is a linearly independent set. We
consider $\Lambda _{j+1}$.We have
$$\Lambda _{j+1}=\{\varphi _{j+1,i}|i=0,1,\cdots,p-1-j-1\}.$$If it
were not a linearly independent set,there would exist some complex numbers
$c_{i}(i=0,1,\cdots,p-1-(j+1))$,at least one of which is nonzero, such
that
$$\sum_{i=0}^{p-1-(j+1)} c_{i}\varphi _{j+1,i}=0.$$
This means 
$$f\sum_{i=0}^{p-1-(j+1)} c_{i}\varphi _{j,i}=0.$$
Now considering $f\varphi _{j,p-1-j}=0$ one can easily convince oneself that
$\sum_{i=0}^{p-1-(j+1)} c_{i}\varphi _{j,i}$ must be a scalar multiple of
$\varphi _{j,p-1-j}$ when Condition 1 is satisfied.
But this is impossible because $\Lambda _{j}$ is
a linearly independent set. The inductive step is thus completed and the lemma
proved.

Theorem 5.1.The set $\{e_{0}^{i}f^{j}\Omega _{0}|i,j=0,1,
\cdots,p-1\}$ is a basis of $V_{\lambda _{1}}\otimes V_{\lambda _{2}}$
if and only if Condition 3 is satisfied.

Proof.Lemma 5.1 implies the "only if" part.The "if" part follows from
Lemma 5.2 since the set   $\{e_{0}^{i}f^{j}\Omega _{0}|i,j=0,1,\cdots,p-1\}$
is none other than $\bigcup_{l=0}^{p-1} \Delta _{l}\bigcup \Lambda _{l}$.

Dually we have the following

Theorem 5.2.The set $\{f_{0}^{i}e^{j}\Phi _{0}|i,j=0,1,
\cdots,p-1\}$ is a basis of $V_{\lambda _{1}}\otimes V_{\lambda _{2}}$
if and only if Condition 2 is satisfied.

We omit the details.

\vspace{6mm}
Acknowledgement. The author thanks Dr. Omar Foda for the hospitality
extended to him during his stay at The University of Melbourne.

\vspace{5mm}
References
\begin{enumerate}
\item V.Chari and A.N.Pressley,Quantum affine algebras,Commun.Math.Phys.
142(1991)261.
\item V.Chari and A.N.Pressley,Quantum affine algebras and their
representations,Canadian Math.Soc.Conf.Proc.16(1995)59.
\item V.Chari and A.N.Pressley,Minimal affinizatios of representations
of quantum groups:the simply-laced case,J.Algebra 184(1996)1.
\item V.Chari and A.N.Pressley,Quantum affine algebras at roots of unity,
q-alg/9609031.
\item J.Beck and V.G.Kac,Finite dimensional representations of quantum affine
algebras at roots of unity,J.Amer.Math.Soc.9(1996)391.
\item V.Chari and A.N.Pressley, A Guide to Quantum Groups,Cambridge University
Press,Cambridge,1994.
\item M.Jimbo,A q-analogue of U(gl(N+1)),Hecke algebra and the Yang-Baxter
equation,Lett.Math.Phys.11(1986)247.
\item X.F.Liu,A new aspect of representations of $U_q(\hat{sl_2})$---generic
case,preprint.
\end{enumerate}
\end{document}